\documentclass[12pt]{article}
\usepackage{psfrag,oldgerm}
\usepackage{amsmath,amsthm,amssymb,bbm,wrapfig}
\usepackage{epsfig,fullpage} 
\newtheorem{theorem}{Theorem}
\newtheorem{corollary}[theorem]{Corollary}

\newtheorem{lemma}[theorem]{Lemma}
\newtheorem{conjecture}[theorem]{Conjecture}
\theoremstyle{remark}

\newtheorem{definition}[theorem]{\bf Definition}

\numberwithin{theorem}{section}

\def\L{\mathfrak{L}}

\newcommand{\intv}{\mathrm{Int}}

\newcommand{\ignore}[1]{}
\newcommand{\NN}{\mathbbm{N}}

\begin{document}
\begin{center}
\textbf{\Large On the Linear Intersection Number of Graphs}\\[2ex]
\textsc{Hauke Klein\\Marian Margraf}\\[2ex]

Mathematisches Seminar,\\
Christian-Albrechts-University of Kiel, Ludewig-Meyn-Str. 4, D-24098
  Kiel, Germany\\
and\\
Institute for Computer Science and Applied Mathematics,\\
  Christian-Albrechts-University of Kiel, Olshausenstr. 40, D-24098
  Kiel, Germany\\[2ex]
klein@math.uni-kiel.de and mma@informatik.uni-kiel.de
\end{center}

\noindent\textbf{Abstract }The celebrated
Erd\"os, Faber and Lov\'asz Conjecture may be
stated as follows: Any linear hypergraph on $v$ points has chromatic index
at most $v$.
We will introduce the linear intersection number of a graph, and use
this number to give an alternative formulation of the Erd\"os, Faber,
Lov\'asz conjecture. Finally, first results about the linear
intersection number will be proved. For example, the definition of the
linear intersection number immediately yields an easy upper bound,
and we determine all graphs for which this bound is sharp.



\section{Introduction}

\begin{definition} A \emph{linear hypergraph} or \emph{partial linear space} is
a pair $\pi=(P,{\cal L})$ 
consisting of a set $P$ of elements called \emph{points} and a set ${\cal L}$ of
distinguished subsets of $P$, called \emph{lines} or
\emph{hyperedges}, satisfying the following
axioms.

\begin{description}
\item[(L1)] Any two distinct points belong to at most one line.
\item[(L2)] Any line has at least two points.
\end{description}
Moreover, a \emph{linear space} is a linear hypergraph in which any two points 
belong to precisely one line. Dually, the linear hypergraph is intersecting if
distinct lines always intersect.
\end{definition}

A \emph{line coloring} $c$ of a linear hypergraph $\pi=(P,{\cal L})$ is a map
$c:{\cal 
L}\longrightarrow \mathfrak{C}$ into some color set $\mathfrak{C}$ such that
any pair of intersecting lines has different colors, i.e. given $l,g\in {\cal
L},\ l\not=g$ then $c(l)\not=c(g)$ if $l\cap g\not=\emptyset.$ The coloring $c$
will be called a $v-$coloring if $|\mathfrak{C}|\leq v.$ Clearly, we are
interested in the minimum cardinality of $\mathfrak{C},$ denoted by
$\chi'(\pi),$ the so-called  \emph{chromatic index} of $\pi.$  

A famous conjecture of Erd\"os, Faber and Lov\'asz can be stated as follows, see
for example \cite{Hin81}.

\begin{conjecture} Every finite
linear hypergraph $\pi=(P,{\cal L})$ with $v=|P|$ points admits a
$v-$coloring of its lines.
\end{conjecture}

The history of results on this conjecture is rather brief. It is known that every
linear hypergraph has a $(\frac{3}{2}\cdot v-3)-$line coloring, 
\cite{CL89}. Moreover, Kahn showed in  \cite{Kahn97} that the conjecture is
asymptotically true, i.e. there is a $(v+o(1))-$line coloring.

The definition of the linear intersection number $v(G)$ for
every graph $G$ will be done in Section 3 and we show that the conjecture of
Erd\"os, Faber and Lov\'asz is 
true if and only if $\chi(G)\leq v(G)$ for all graphs $G$, where $\chi(G)$
denotes the chromatic number of $G$. 
Moreover, Section 3 contains general auxiliary results about the
linear intersection number, in particular we will determine all
graphs with maximal linear intersection number.
Finally, Section 4 contains theorems concerning
lower bounds for the linear intersection number.

\section{Notations}

Let $\pi=(P,{\cal L})$  be a linear hypergraph, $v=|P|$ and $b=|{\cal
L}|$. By ${\cal
L}_p=\{l\in{\cal L};p\in l\}$ we denote the \emph{line pencil} of $p\in P$ and
by $r_p=|{\cal L}_p|$ the \emph{degree} of $p.$ Dually, $k_l=|l|$ denotes the
cardinality of a line $l\in {\cal L}.$ The hypergraph $\pi$ is called
$r-$uniform if $k_l=r$ for all $l\in {\cal L}.$ A \emph{point clique}  is a
subset $C$ of 
points such that any two points of $C$ are joined by a line. Dually, a
\emph{line clique} is a subset ${\cal C}$ of lines such that any two lines of
${\cal C}$ intersect. The \emph{clique number} $\omega(\pi)$ is the
maximum cardinality of a point clique of $\pi$, while the \emph{clique
index} $\omega'(\pi)$ is the
maximum cardinality of a line clique.

An important theorem which we will use frequently is the so-called Fundamental
Theorem of finite linear spaces, see \cite{BB93} Theorem 1.5.5.

\begin{theorem} Let $\pi$ be a finite linear space. Then $b\geq v.$
Moreover, equality holds if and only if $\pi$ is a projective plane or a
near pencil.
\end{theorem}

Applying this to the dual hypergraph $\pi^*$, we obtain that $v$ is an upper bound
of the clique index (see for example \cite{Kahn97}).

\begin{corollary} \label{clique2}If $\pi$ is intersecting then $b\leq v,$ i.e. the
conjecture is true in this case. 
\end{corollary}

\section{The linear intersection number}
Given any linear hypergraph $\pi=(P,{\cal L})$, we define the 
\emph{intersection graph} $G_{\pi}$ of $\pi$ to be the graph
$G_\pi=({\cal L},E)$ whose edges are the pairs of intersecting
lines in $\pi$, i.e.
$E=\{\{l,g\}; l\cap g\not=\emptyset\}$. In other words, two lines
are joined in $G_{\pi}$ if and only if they have a common point
in the hypergraph $\pi$.
Obviously, a map $c:V\mapsto \mathfrak{C}$ is a
vertex coloring of $G_\pi$ if and only if $c$ is a line coloring of $\pi$.
Hence, to determine the chromatic index of a hypergraph $\pi$, it is
enough to look at its intersection graph $G_{\pi}$.
But unfortunately, it is not possible to recover the original
hypergraph by looking only at its intersection graph. In fact,
we cannot even see the number $v$ of points of $\pi$ given only
$G_{\pi}$.
In order to solve this problem, we are going to
define the linear intersection
number of an arbitrary graph.

Let $G=(V,E)$ be a graph. Given any vertex $x\in V$, let $E_x:=\{e\in E;
x\in E\}$ be the set of all edges incident with $x$.
Then the dual space
\[
G^*=(E,\{E_x;x\in V\})
\]
of $G$
forms a hypergraph with $G_{G^*}=G$.
Moreover, this hypergraph is almost a linear hypergraph.
Any two lines of $G^*$
intersect in at most one point, however any vertex $x\in V$ of
degree one yields a line $E_x$ of $G^*$ with only one point,
and such lines are not allowed in a linear hypergraph.
But this is not a real problem, just add a new point to any
line $E_x$ with $r_x=1$, and any vertex $x\in V$ with $r_x=0$
needs two points for its line. The linear hypergraph 
$\overline{G}^*$ obtained by this construction has our original
graph as its intersection graph, i.e. $G_{\overline{G}^*}=G$.

In particular, each graph is the intersection graph of a
linear hypergraph and the set
\[
\mathfrak{P}_G:= \{ v\in \NN;
\pi \mbox{ is a linear hypergraph on $v$ points such that
}G_\pi=G\}
\]
is not empty. As usual, the equality $G_{\pi}=G$ in this formula
actually means isomorphism.
We will define the linear intersection number of $G$ to be the
smallest number of points of a linear hypergraph realizing $G$
as its intersection graph.

\begin{definition} For every  graph $G=(V,E)$ let 
$v(G):=\min\mathfrak{P}_G$ be
the  \emph{linear intersection number}.
\end{definition}

Now we have a nice description of the Erd\"os, Faber and Lov\'asz conjecture in
terms of the chromatic number of graphs.
\begin{corollary} The Erd\"os, Faber and Lov\'asz conjecture is true if and
only if
$\chi(G)\leq v(G)$ for any graph $G$.
\end{corollary}

{\it Proof:} Assume the conjecture is true. Let $G$ be a 
graph. Then for every linear hypergraph $\pi=(P,{\cal L})$ with
$G_\pi=G$ we have $|P|\geq \chi'(\pi)=\chi(G_\pi)=\chi(G),$
hence $v(G)\geq \chi(G).$ 

Conversely, let $G$ be a  graph with $\chi(G)\leq v(G).$ Then we obtain
for every linear hypergraph $\pi=(P,{\cal L})$ with $G_\pi=G$ that
$|P|\geq v(G)\geq \chi(G)=\chi'(\pi).$\qed

Since we explicitly constructed a linear hypergraph with given
intersection graph, we get a simple upper bound on $v(G)$. In
order to formulate this result, define for any graph $G=(V,E)$
\begin{eqnarray*}
L(G)&:=&\{x\in V;r_x=1\}\quad\text{the set of leafs of $G$},\\
l(G)&:=&|\{x\in V;r_x=1\}|\quad\text{the number of leafs of $G$},\\
I(G)&:=&\{x\in V;r_x=0\}\quad\text{the set of isolated vertices of $G$ and}\\
e(G)&:=&|\{x\in V;r_x=0\}|\quad\text{the number of isolated vertices of $G$.}
\end{eqnarray*}
Since $\overline{G}^*$ forms a linear hypergraph satisfying
$G_{\overline{G}^*}=G$  we already know
\begin{corollary}\label{edgebound}
Given any graph $G$, we have $v(G)\leq |E|+l(G)+2e(G)$.
\end{corollary}

The inequality in this corollary may actually become an equality,
and we will determine all graphs $G$ with $v(G)=|E|+l(G)+2e(G)$
later in this section.
The linear intersection number is additive, i.e. given two graphs
$G_1,G_2$ and a realisation of $G_1+G_2$ as intersection graph of
a linear hypergraph, no line of $G_1$ can intersect a line of $G_2$,
hence the linear hypergraph is a sum of two linear hypergraphs
realizing $G_1$ and $G_2$, respectively. This observation proves
\[
v(G_1+G_2)=v(G_1)+v(G_2).
\]
As it turns out, it is possible to describe the linear intersection
number without even mentioning linear hypergraphs.

\begin{lemma}\label{cliqchar}
Let $G=(V,E)$ be a graph. Then
\[
v(G)=\min\left\{r\ge0\left|\begin{array}{l}
\text{There exist cliques $C_1,\ldots,C_r\subseteq V$ of $G$ such that}\\
\text{any edge of $G$ is in exactly one $C_i$ and each vertex}\\
\text{of $G$ belongs to at least two of these cliques}
\end{array}\right.\right\}.
\]
\end{lemma}

{\it Proof.} Write $v=v(G)$ and let $\pi=(P,V)$ be a linear hypergraph
with $G_{\pi}=G$ and $|P|=v$. Given any point $x\in P$ of $\pi$, the
set $C_x:=\{a\in V;x\in a\}$ of all lines passing through $x$ is a clique
in $G=G_{\pi}$. An edge $e$ in $G$ is a pair $(a,b)$ of intersecting
lines $a,b$ of $\pi$, but these two lines intersect in exactly one point
$x\in P$, hence there is exactly one $C_x$ containing $a$ and $b$.
Finally, each line of $\pi$ contains at least two points, hence each
vertex of $G$ is contained in at least two $C_x$. 

Conversely, assume we got a collection $C_1,\ldots,C_r$ of cliques
of $G$, as above. Build a linear hypergraph by using
$P:=\{C_1,\ldots,C_r\}$ as its set of points. Given any vertex
$a\in V$, construct a set of points corresponding to $a$ by
$\overline{a}=\{C_i;a\in C_i\}$. Then, $(P,\{\overline{a};a\in V\})$
is a linear hypergraph with intersection graph $G$ and $|P|=r$.\qed

It is quite natural to remove the condition that each vertex belongs
to at least two of the cliques, i.e. to consider arbitrary partitions
of the set of edges of $G$ into cliques. In this context, there is
no need for trivial cliques, i.e. cliques with at most one element,
and we define
\[
{\cal C}(G):=\{C\subseteq V;\text{$C$ is a clique with $|C|>1$}\}
\]
to be the set of all non-trivial cliques of a graph $G=(V,E)$.
Now, we will define the reduced linear intersection number of
a graph.

\begin{definition}Let $G=(V,E)$ be a graph. The
\emph{reduced linear intersection number of $G$} is defined to be
\[
\overline{v}(G):=\min\left\{r\ge0\left|\begin{array}{l}
\text{There exist $C_1,\ldots,C_r\in{\cal C}(G)$ such that}\\
\text{each edge of $G$ is in exactly one $C_i$}
\end{array}\right.\right\}.
\]
\end{definition}

Obviously, $\overline{v}(G)\le v(G)$ for any graph $G$, but
the reduced linear intersection number may be quite
different from $v(G)$, for example $\overline{v}(G)=1$ if $G$
is a complete graph. However, in many cases the two numbers will
coincide. Differences between $v(G)$ and $\overline{v}(G)$ are
caused by the existence of some kind of interior vertices of the
graph $G$. We are going to define a whole bunch of variants of
the concept of an interior vertex of $G$.

\begin{definition}Let $G=(V,E)$ be a graph, and let $x\in V$ be a vertex
of $G$. The neighborhood of $x$ will be denoted by $G_x:=\{y\in V;
\{x,y\}\in E\}$.\begin{enumerate}
\item The vertex $x$ is an \emph{interior vertex} of $G$ if $G_x$ is a clique.
Let $\intv(G)$ be the set of all interior vertices of $G$.
\item The set of \emph{strongly interior vertices} of $G$ is defined
to be
\[
\intv_s(G):=\{x\in\intv(G);|G_x|>1\},
\]
i.e. $\intv_s(G)=\intv(G)\backslash(L(G)\cup I(G))$.
\item An interior vertex $x$ of $G$ is an \emph{extremal interior vertex}
if there exists a collection $C_1,\ldots,C_r$ as in Lemma \ref{cliqchar}
with $r=v(G)$ and $G_x\cup\{x\}=C_i$ for some $1\le i\le r$. Of course,
$G_x\cup\{x\}=C_i$ for some $i$ if and only if $\{x\}=C_j$ for some $j$.
The set of
all extremal interior vertices of $G$ will be denoted by $\intv_e(G)$.
Note that leafs and isolated vertices of $G$ are also
extremal interior vertices of $G$.
\item Finally, define
\[
\intv_{es}(G):=\intv_s(G)\cap\intv_e(G).
\]
\end{enumerate}\end{definition}

Of course, it may be difficult to recognize the extremal interior vertices
in the set of all interior vertices of a graph. We are going to discuss
all these concepts in a complete graph. As usual, the maximal size of a
clique in a graph $G$ will be denoted by $\omega(G)$.
Each clique in the intersection graph of a linear hypergraph
corresponds to a set of pairwise intersecting lines in this linear
hypergraph, hence
Corollary \ref{clique2} immediately implies
a lower bound for the linear intersection number.

\begin{theorem} \label{th44} For every graph $G$ we have
$\omega(G)\leq v(G).$ 
\end{theorem}
\qed 

Now, it is easy to discuss the complete graph.

\begin{theorem}\label{kntheorem}
Let $G=K_n$ be a complete graph on $n\ge3$ vertices.
Then $v(K_n)=n$, $\overline{v}(K_n)=1$ and $\intv_e(K_n)=\emptyset$.
Moreover, $v(K_1)=2$ and $v(K_2)=3$.\end{theorem}

\textit{Proof.} Theorem \ref{th44} implies $n=\omega(K_n)\le v(K_n)$.
Choose a vertex $a\in V$ where $K_n=(V,E)$, and consider the following
collection of cliques of $K_n$
\[
V\backslash\{a\},(\{a,x\})_{x\in V\backslash\{a\}}.
\]
Lemma \ref{cliqchar} yields $v(K_n)\le n$. 
Finally, assume $C_1,\ldots,C_s$ is a collection of cliques of $K_n$
such that each edge is contained in exactly one $C_i$ and each vertex
is contained in at least two $C_i$. If $C_i=V$ for some $i$, then
$|C_j|\le1$ for $j\not=i$, and in particular $s\ge n+1$. Consequently,
there are no extremal interior vertices in $K_n$.\qed

Geometrically, the linear hypergraph on $n$ points realizing $K_n$
is a so-called near-pencil.
\begin{center}
\epsfig{file=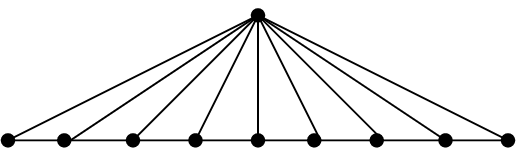}

{\small The near-pencil on ten points}
\end{center}
The result $\intv_e(K_n)=\emptyset$ just states that the set of all
vertices of $K_n$ cannot occur in a minimal partition of the set of
all edges of $K_n$ into cliques. It is tempting to ask which subsets
of $K_n$ actually occur in such a decomposition. Of course, an answer
depends only on the size $r$ of such a set. Our construction gives
a positive answer if $r=2$ or $r=n-1$ and negative answers if $r=n$ or
$r=1$. Other values of $r$ correspond to projective planes with $n$
points, i.e. $r\not=2,n-1$ is possible if and only if $n=k^2+k+1$
and there exists a projective plane of order $k$. In this case, $r$
is the size of a line pencil, i.e. we have $r=k+1$.

We will prove two simple results relating $v(G)$ and $\overline{v}(G)$.

\begin{lemma}\label{restrict}
Let $G=(V,E)$ be a graph.\begin{enumerate}
\item If $\intv_{es}(G)=\emptyset$, then $v(G)\ge\overline{v}(G)+l(G)+2e(G)$.
\item If $\intv_s(G)=\emptyset$, then $v(G)=\overline{v}(G)+l(G)+2e(G)$.
\end{enumerate}
In particular, $\overline{v}(G)=v(G)$ if $G_x$ is not a clique
for any vertex $x$ of $G$.\end{lemma}

\textit{Proof.} First, assume there are no extremal, strongly
interior vertices of $G$.
Let $C_1,\ldots,C_{v(G)}$ be a sequence of cliques
of $G$ such that each edge is contained in exactly one of these
cliques, and each vertex is contained in at least two of them.
We may assume that exactly the cliques $C_1,\ldots,C_s$, $0\le s\le r$
are non-trivial, and in particular $\overline{v}(G)\le s$. Each
of the remaining cliques $C_{s+1},\ldots,C_{v(G)}$ consists of
exactly one point. Let $s+1\le i\le r$ and $C_i=\{x\}$. By
minimality, $x$ is contained in at most one of the cliques
$C_1,\ldots,C_s$.

We claim that $x\in L(G)\cup I(G)$. In fact, if $x\notin C_1,\ldots C_s$,
then there is no edge incident with $x$ and $x$ is an isolated vertex.
Otherwise $x\in C_j$ for exactly one $1\le j\le s$, and each edge incident
with $x$ is in $C_j$, i.e. $G_x\cup\{x\}\subseteq C_j$. Hence $C_j=
G_x\cup\{x\}$ and $x\in\intv_e(G)$. But $\intv_{es}(G)=\emptyset$, i.e.
$|G_x|=1$ and $x$ is a leaf of $G$.

If $x$ is an isolated vertex, then $x\notin C_1,\ldots,C_s$, so
there is another $s+1\le j\le v(G)$ with $j\not=i$ and $C_j=\{x\}$.
This proves $v(G)-s\ge l(G)+2e(G)$, i.e. $v(G)\ge s+l(G)+2e(G)\ge
\overline{v}(G)+l(G)+2e(G)$.

This argument proves the first assertion. Now, assume $\intv_s(G)=
\emptyset$. Let $r:=\overline{v}(G)$, and let $C_1,\ldots,C_r$
be a collection of non-trivial cliques of $G$ such that each edge
is contained in exactly one of these cliques. Given any vertex
$x\in V\backslash(L(G)\cup I(G))$ which is neither a leaf nor an isolated
vertex, we choose an edge $e\in E_x$ and an $i$ such that $e$ is an
edge of $C_i$. Since $G_x$ is not a non-trivial clique and $r_x>1$,
we know that $E_x\not\subseteq C_i$, i.e. there is a $j\not=i$
containing another edge $f\in E_x$. In particular, the vertex $x$
is in at least two of $C_1,\ldots,C_r$. Consequently,
\[
(C_i)_{1\le i\le r},\,(\{x\})_{x\in L(G)},\,(\{x\},\{x\})_{x\in I(G)}
\]
is a collection of cliques of $G$ as in Lemma \ref{cliqchar}, and we
have $v(G)\le r+|L(G)|+2|I(G)|=\overline{v}(G)+l(G)+2e(G)$.\qed

\begin{lemma}Let $G=(V,E)$ be a graph on $n=|V|$ vertices.
Then $v(G)\le\overline{v}(G)+n+e(G)$ and equality holds if
and only if $G$ is a sum of isolated vertices and edges.
\end{lemma}

\textit{Proof.} Let $r=\overline{v}(G)$ be the reduced linear intersection
number of $G$, and choose non-trivial cliques $C_1,\ldots,C_r\subseteq V$
such that each edge of $G$ is in exactly one $C_i$.
Let $S$ be the set of vertices of $G$ contained in exactly
one $C_i$. By adding one trivial clique $\{x\}$ for each $x\in S$ and
two trivial cliques $\{x\}$ for each $x\in I(G)$, we obtain $v(G)\le r+|S|+
2|I(G)|\le r+n+e(G)$.
\par Now assume, $v(G)=r+n+e(G)$, in particular $|S|+2|I(G)|=n+|I(G)|$, i.e.
$S=V\backslash I(G)$ and $C_1,\ldots,C_r$ is a partition of $V\backslash I(G)$
into cliques. Since each edge of $G$ is contained in some $C_i$,
there are no edges between $C_i$ and $C_j$ if $i\not=j$, hence $G$
is a sum of complete graphs.
\par Finally assume $G=K_1+\ldots+K_s+L_1+\ldots+L_t+I(G)$ for complete graphs
$K_1,\ldots,K_s$, $L_1,\ldots,L_t$ with $|K_i|\ge3$, $|L_j|=2$.
Then $\overline{v}(G)=s+t$, and $v(G)=|K_1|+\ldots+|K_s|+3t+2|E|=
n+t+e(G)$ by Theorem \ref{kntheorem},
consequently $v(G)=\overline{v}(G)+n+e(G)$
if and only if $s=0$.\qed

Let $G=(V,E)$ be a graph. We will build a new graph $C_G$ by
taking the set ${\cal C}(G)$ of all non-trivial cliques of
$G$ as its set of vertices. Two non-trivial cliques $A,B\subseteq V$
of $G$ will be joined by an edge in $C_G$ if and only if they
have a common edge, i.e. $|A\cap B|\ge2$. The graph $C_G$ will
be called the \emph{clique graph of $G$}. The stability number
$\alpha(H)$ of a graph $H$ is defined to be the maximal size
of an independent set of vertices of $H$.

\begin{corollary}Let $G=(V,E)$ be a graph. Then
\[
\overline{v}(G)=\min\{|U|;
\text{$U\subseteq{\cal C}(G)$ is a maximal independent set in $C_G$}\}
\le\alpha(C_G).
\]
\end{corollary}

\textit{Proof.} Let $r=\overline{v}(G)$ and choose $C_1,\ldots,C_r\in
{\cal C}(G)$ as in the definition of $\overline{v}(G)$. Then
$\{C_1,\ldots,C_r\}$ is a maximal independent subset of $C_G$ since
each edge of $G$ is contained in one of the $C_i$.
Conversely, let $\{C_1,\ldots,C_s\}$ be a maximal
independent set of vertices of ${\cal C}(G)$.
If $e$ is an edge of $G$ not in any of the $C_i$,
then the two vertices incident with $e$ form another non-trivial
clique which could be added to $C_1,\ldots,C_s$ contradicting the
maximality assumed of $\{C_1,\ldots,C_s\}$.
Hence, each edge of $G$ is in exactly one $C_i$,
and $r\le s$.\qed

It seems to be important to study the behaviour of the intersection
number of a graph with respect to natural constructions on graphs.
Let $G=(V,E)$ be a graph. Given two vertices $a,b\in V$ not incident
with a common edge, i.e. $(a,b)\notin E$, let $G/ab$ be the graph
obtained from $G$ by collapsing the set $\{a,b\}$ into a single vertex.

\begin{lemma}\label{collaps}
Let $G=(V,E)$ be a graph, and $a,b\in V$ be two vertices of $G$.
\begin{enumerate}
\item If $d(a,b)\ge3$, then $v(G/ab)\le v(G)$.
\item If $d(a,b)\ge4$ and $a,b\notin\intv(G)$, then $v(G/ab)=v(G)$.
\end{enumerate}\end{lemma}

\textit{Proof.} Denote the set of all vertices of $G/ab$
by $V^{\prime}=V\backslash\{a,b\}\cup\{\omega\}$.
We begin with the first result. Let $\pi=(P,V)$ be a linear hypergraph
with intersection graph $G$ and $|P|=v(G)$. Replace the lines $a,b$ of
$\pi$ by a new line $\omega=a\cup b$. Since $d(a,b)\ge3$, there is no
line in $\pi$ intersecting $a$ and $b$, hence $\pi^{\prime}=(P,V^{\prime})$ is a
linear hypergraph with $G_{\pi^{\prime}}=G/ab$, in particular $v(G/ab)
\le|P|=v(G)$.

Now, assume $d(a,b)\ge4$ and $a,b\notin\intv(G)$. Let $\pi=(P,V^{\prime})$
be a linear hypergraph with intersection graph $G/ab$ and $|P|=v(G/ab)$.
The neighborhood of $\omega$ in $G/ab$ is the union $G_a\cup G_b$, and
$d(a,b)\ge3$ implies $G_a\cap G_b=\emptyset$. Define a new hypergraph
$\pi^{\prime}=(P,V)$ by replacing the line $\omega$ with the two new lines
$a=\{l\cap\omega;l\in G_a\}$ and $b=\{l\cap\omega;l\in G_b\}$. We claim
that $\pi^{\prime}$ is a linear hypergraph.

Given points $p,q\in P$ and lines $m,n\in V$ with $p\not=q$, $p,q\in l,m$,
we shall show that $m=n$. If $m,n\not=a,b$, we immediately have $m=n$ since
$\pi$ is a linear hypergraph. Otherwise, we may assume that $m=a$.
In particular $p,q\in\omega$, hence $n\notin V\backslash\{a,b\}$, i.e.
$n=a$ or $n=b$. If $n=b$, then $p=l\cap\omega=l^{\prime}\cap\omega$ for some
$l\in G_b$, $l^{\prime}\in G_a$, hence $d(a,b)\le3$; a contradiction.
This contradiction implies $n=a=m$. It remains to prove that $a$ and $b$
are each incident with at least two points. Since $a$ is not an isolated
vertex of $G$, there is a vertex $l\in G_a$, hence $l\cap\omega\in a$. If
$l\cap\omega$ is the only point of $a$, then $l^{\prime}\cap\omega=
l\cap\omega$ and $(l,l^{\prime})\in E$ for each $l^{\prime}\in G_a$,
$l^{\prime}\not=l$, i.e. $G_a$ is a clique and $a$ is an interior vertex
of $G$; a contradiction. Similarly, the line $b$ has at least $2$ point.

Finally, $G_{\pi^{\prime}}=G$ implies $v(G)\le|P|=v(G/ab)$.\qed

We will need the special case $d(a,b)=\infty$ of this lemma, i.e.
the vertices are in distinct connected components of $G$. In this
special case, we have a more complete result.
It is usefull to introduce a small notation. Assume, $G_i=(V_i,E_i)$
are two graphs with a common vertex $V_1\cap V_2=\{a\}$. The join of
$G_1$ and $G_2$ at $a$ is the graph $G_1\vee_a G_2:=(V_1\cup V_2,E_1\cup E_2)$,
i.e. the two graphs $G_1$ and $G_2$ are glued along the vertex $a$.
Since each non-trivial clique in $G_1\vee_aG_2$ is contained in
exactly one of the two graphs $G_1$ and $G_2$, we immediately see
\[
C_{G_1\vee_aG_2}=C_{G_1}+C_{G_2}\quad\text{and}\quad
\overline{v}(G_1\vee_aG_2)=\overline{v}(G_1)+\overline{v}(G_2).
\]
However, the intersection number of $G_1\vee_aG_2$ is slightly more
complicated.

\begin{lemma}\label{gluing}
Let $G_i=(V_i,E_i)$, $i=1,2$ be two graphs with a common vertex
$V_1\cap V_2=\{a\}$. Write $t:=|\{i\in\{1,2\};a\in\intv_e(G_i)\}|$.
Then $v(G_1\vee_aG_2)=v(G_1)+v(G_2)-2$ if $a$ is either an isolated
vertex of $G_1$ or an isolated vertex of $G_2$. Otherwise, $v(G_1\vee_aG_2)=
v(G_1)+v(G_2)-t$ and $a\notin\intv_e(G)$.
\end{lemma}
\textit{Proof.} If $a$ is an isolated vertex of $G_2$, write $G_2=a+G_3$,
and compute $v(G_1\vee_aG_2)=v(G_1+G_3)=v(G_1)+v(G_3)=v(G_1)+v(G_2)-2$.
Hence we may assume that $a\notin I(G_1)\cup I(G_2)$.
For $i=1,2$, set $v_i:=v(G_i)$ and choose cliques
$C_1^i,\ldots,C^i_{v_i}$ as in Lemma \ref{cliqchar}. Moreover,
assume $C^i_j=\{a\}$ for some $j$ if $a\in\intv_e(G_i)$. Then
form the collection of all $C^i_k$ omitting each occurrence of
$\{a\}$. These are $v_1+v_2-t$ cliques in $G_1\vee_aG_2$, and
Lemma \ref{cliqchar} implies $v(G_1\vee_aG_2)\le v_1+v_2-t$.

Conversely, let $r=v(G_1\vee_aG_2)$ and let $C_1,\ldots,C_r$ be cliques
in $G_1\vee_aG_2$ such that each edge of $G_1\vee_aG_2$ is in exactly
one $C_i$ and each vertex of $G_1\vee_aG_2$ is in at least two of the
$C_i$. For $i=1,2$ set $t_i:=1$ if $a$ is an extremal interior vertex
of $G_i$, and $t_i:=0$ otherwise.
Since $a$ is not an isolated vertex of $G_1$ or $G_2$ we have
$C_i\not=\{a\}$ for each $1\le i\le r$. Assume that $C_1,\ldots,C_s
\subseteq V_1$ and $C_{s+1},\ldots,C_r\subseteq V_2$.

Each edge of $G_1$ is contained in exactly one of the $C_1,\ldots,C_s$,
and each vertex $x\in V_1\backslash\{a\}$ is contained in at least two
of the $C_1,\ldots,C_s$. In particular, $s\ge v_i-1$. If $s=v_1-1$,
then $a\in\intv_e(G_1)$ and $t_1=1$, i.e. we know $s\ge v_1-t_1$ in
any case. Similarly $r-s\ge v_2-t_2$, i.e. $r\ge v_1+v_2-(t_1+t_2)=
v_1+v_2-t$.\qed

We will use this lemma to discuss the effect of removing a clique
from a given graph.

\begin{lemma}\label{cliqremove}
Let $G=(V,E)$ be a graph and let $C\subseteq V$ be a clique of
$G$ with $n=|C|\ge3$. Let $G^-:=G-C=(V^-,C^-)$ be the residual graph, i.e.
\begin{eqnarray*}
E^-&:=&\{e\in E;e\not\subseteq C\},\\
V^-&:=&\bigcup E^-\cup(V\backslash C).
\end{eqnarray*}
Then $v(G)\le v(G^-)+n$ and $v(G)=v(G^-)+n$ if and only if
either $C$ is a connected component of $G$ or
$V^-\cap C=\{a\}$ for a vertex $a\notin\intv_e(G^-)$.\end{lemma}

\textit{Proof.} Write $r:=|V^-\cap C|$. If $r=0$, then there are no
edges between $C$ and $V\backslash C$, i.e. $C$ is a connected
component of $G$, $G=C+G^-$ and $v(G)=v(C)+v(G^-)=n+v(G^-)$ by
Theorem \ref{kntheorem}.

Now, assume that $r\ge2$. Write $s=v(G^-)$ and choose cliques
$C_1,\ldots,C_s$ in $G^-$, as usual. Looking at $C_1,\ldots,C_s,
C,(\{x\})_{x\in C\backslash V^-}$, Lemma \ref{cliqchar} yields
$v(G)\le s+1+n-r\le s+n-1<v(G^-)+n$.

Finally, assume $r=1$, i.e. $G=C\vee_aG^-$. Then $v(G)=v(G^-)+n-t$
by Lemma \ref{gluing} and Theorem \ref{kntheorem} where $t=1$ if
$a\in\intv_e(G^-)$ and $t=0$ if $a\notin\intv_e(G^-)$.\qed

Now, we are in a position to discuss the case of equality in
Corollary \ref{edgebound}. First, we shall describe the class
of graphs realizing equality.

\begin{definition}A graph $G$ is called almost triangle-free if
is obtained using the following construction. Begin with a triangle-free
graph $G_0=(V_0,E_0)$. Glue an arbitrary number of triangles $D_i$ at each
vertex $a\in V_0$, but these triangles shall not intersect each other
except in the one vertex glued to $G_0$. Vertices in the triangles $D_i$
but not in $G_0$ are called extremal vertices of $G$ unless the triangle
$D_i$ is the only triangle glued to an isolated vertex of $G_0$.
\end{definition}

For example, the graph

\begin{center}
\epsfig{file=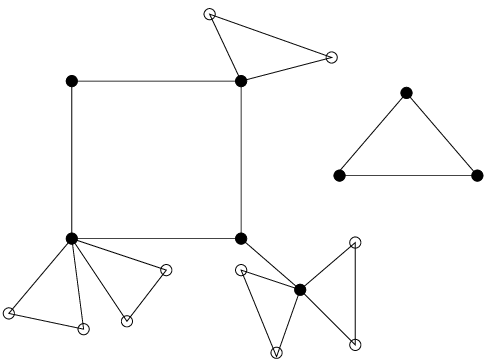}
\end{center}

is almost triangle-free, and the vertices marked
with an open circle are the extremal vertices of this graph.
However, the graph

\begin{center}
\epsfig{file=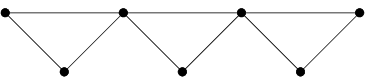}
\end{center}

is not almost triangle-free.

\begin{lemma}\label{almosttf}Let $G$ be an almost triangle-free graph with
$m$ edges. Then $v(G)=m+l(G)+2e(G)$ and a vertex $x$ of $G$ is an extremal,
strongly interior vertex of $G$ if and only if $x$ is an exterior vertex.
\end{lemma}

\textit{Proof.} First, assume $G=(V,E)$ is a triangle-free graph with
$m=|E|$ edges. Then each non-trivial clique of $G$ has exactly two
vertices, hence $\overline{v}(G)=m$.
Moreover, $\intv(G)=L(G)\cup I(G)$
and $\intv_{es}(G)=\intv_s(G)=\emptyset$. In particular,
$v(G)=m+l(G)+2e(G)$ by Lemma \ref{restrict}.

Now, assume $G$ is an almost triangle-free graph with $m$ edges and our lemma
holds for $G$. Let $a$ be a vertex of the triangle-free graph $G$
was build from, and let $D$ be a new triangle intersecting $G$ in
the vertex $a$. We shall show that the lemma remains to be true for
the new graph $G^{\prime}=G\vee_aD$ with $m^{\prime}=m+3$ edges.
We have to distinguish three distinct cases.

If $a$ is an isolated vertex of $G$, then $v(G^{\prime})=v(G)+1=m+1+l(G)+
e(G)=m^{\prime}+l(G^{\prime})+2e(G^{\prime})$. Moreover, $D$ is the only
triangle glued to $a$ in $G^{\prime}$, hence $\intv_{es}(G^{\prime})=
\intv_{es}(G)$ consists exactly of the extremal vertices of $G$, and these
are exactly the extremal vertices of $G^{\prime}$. If $a$ is a leaf of
$G$, then $v(G^{\prime})=v(G)+2=m+2+l(G)+2e(G)=m^{\prime}+l(G^{\prime})+
2e(G^{\prime})$, and the extremal, strongly
interior vertices of $G^{\prime}$ are
exactly the extremal interior vertices of $G$ and the two new
vertices of $D$, and these are exactly the exterior vertices of $G^{\prime}$.

Finally, assume $a\notin L(G)\cup I(G)$. Then $a\notin\intv_e(G)$ by
our hypothesis about $G$, and $v(G^{\prime})=v(G)+3=m+3+l(G)+2e(G)
=m^{\prime}+l(G^{\prime})+2e(G^{\prime})$. Moreover, the strongly, extremal
interior vertices of $G^{\prime}$ are exactly the strongly, extremal interior
vertices of $G$ and the two new vertices of $D$, i.e. the exterior vertices
of $G^{\prime}$.\qed

\begin{theorem}Let $G$ be a graph with $m$ edges. Then
$v(G)=m+l(G)+2e(G)$ if and only if $G$ is almost triangle-free.
\end{theorem}

\textit{Proof.} Let $G$ be a graph on $n$ vertices with $m$ edges,
and assume that $v(G)=m+l(G)+2e(G)$. We shall prove that $G$ is almost
triangle-free. By induction, assume this result
to be true for all graphs with at most $n-1$ vertices. If $G$ is
even triangle-free, we are immediately done. Now, assume that $G$
contains a triangle, i.e. $\omega=\omega(G)\ge3$. Choose a clique
$C$ of $G$ with $|C|=\omega$, and let $G^-=G-C$ be the residual
graph. By Lemma \ref{cliqremove} and Corollary \ref{edgebound},
we know that
\begin{multline*}
m+l(G)+2e(G)=v(G)\le v(G^-)+\omega\le m-\frac{\omega(\omega-1)}{2}+
l(G^-)+2e(G^-)+\omega\\
=m-\frac{\omega(\omega-3)}{2}+l(G)+2e(G).
\end{multline*}
This implies $\omega=3$, $v(G^-)=m-\omega(\omega-1)/2+l(G^-)+2e(G^-)$
and $v(G)=v(G^-)+\omega$. Hence, $G^-$ is almost triangle-free by
induction and either $C$ is a connected component of $G$ or $C$
intersects $G^-$ in a single vertex $a\notin\intv_e(G^-)$. If $C$
is a connected component of $G$, the graph $G=G^-+K_3$ is almost
triangle-free. In the other case, assume the graph $G^-$ was build
from a triangle-free graph $G_0$.
Since $a\notin\intv_{es}(G^-)$, the vertex $a$
is not an exterior vertex of $G^-$. If $a$ is a vertex of $G_0$,
the graph $G$ is almost triangle-free. Otherwise $a$ belongs to
a triangle $D$ of $G^-$ glued as the only triangle to an isolated
vertex of $G_0$. Replacing this isolated vertex by $a$ shows that
$G$ is again almost triangle-free.\qed

In particular, we know the linear intersection number $v(G)$ if
$G$ is a triangle-free graph. This class of graphs includes bipartite
graphs, trees and cycles, so we get a number of immediate corollaries.

\begin{corollary}If $G=(V,E)$ is a bipartite graph with
$m=|E|$ edges, then $v(G)=m+l(G)+2e(G)$.\end{corollary}

\begin{corollary}If $G=(V,E)$ is a tree with $n=|V|$ vertices,
then $v(G)=n+l(G)-1$.\end{corollary}

\begin{corollary} \label{cycle}
For each cycle $C_n,\ n\geq 3$,  we have $v(C_n)=n.$
\end{corollary}

\section{Lower bounds for the linear intersection number}

It is well-known (and very easy to see) that the clique number
is a lower bound for the chromatic number of $G$, and we already
proved that $\omega(G)\le v(G)$.
By Theorem \ref{th44} it follows that our conjecture is true for graphs $G$ 
satisfying $\omega(G)=\chi(G),$ in particular for perfect graphs.
We are going to derive three other results on lower bounds for
the linear intersection number $v(G)$.

First of all, if $\pi$ is a linear hypergraph with $b$ lines on
$v$ points, then $b\le v(v-1)/2$, i.e. $v\ge(1+\sqrt{1+8b})/2$.
Consequently, if $G$ is a graph with $n$ vertices then
\[
v(G)\ge\frac{1+\sqrt{1+8n}}{2},
\]
is a trivial lower bound on $v(G)$. Of course, equality holds if
and only if $G$ is the intersection graph of a complete graph.

A somewhat more interesting bound is given by a theorem of
Seymour on matchings in a linear hypergraph \cite{seymour}.
Here, a matching of a hypergraph is a collection of pairwise
disjoint lines in the hypergraph.

\begin{theorem}If $G$ is a graph with $n$ nodes, then
$\alpha(G)v(G)\ge n$.\end{theorem}

\textit{Proof. }Let $\pi$ be a linear hypergraph on $v(G)$ points
with intersection graph $G=G_{\pi}$. Seymours theorem \cite{seymour}
states the existence of a matching $S$ of $\pi$ consisting of at
least $|S|\ge n/v(G)$ lines. But a matching of $\pi$ is just an
independent set in $G=G_{\pi}$, hence $\alpha(G)\ge n/v(G)$.\qed

Let $\pi=(P,V)$ be a linear hypergraph with intersection graph
$G=(V,E)$. Let $a\in V$ be a line of $\pi$. For any point $x\in a$,
the pencil $V_x=\{b\in V;x\in b\}$ is a clique in $G$ containing $a$.
Moreover, $V_x\cap V_y=\{a\}$ for any two points $x,y\in a$, $x\not=y$.
Hence the sets $V_x\backslash\{a\}$ form a partition of the subgraph
$G_a=\{b\in V;(a,b)\in E\}$ into cliques. Hence, if $\theta(H)$
denotes the minimal number of cliques necessary to partition the
set of vertices of a graph $H$, we have
\[
|a|\ge k(a):=\min\{\theta(G_a),2\}.
\]
The number $k(a)$ is something like the minimal number of points on $a$
in any realization of $G$ as an intersection graph of a linear hypergraph.
Whether there actually exists such a realization $\pi=(P,{\cal L})$
of $G$ with $v(G)=|P|$ and exactly $k(a)$ points on $a$ is another
question. However,
if $U\subseteq V$ is an independent subset of $G$, then the lines of $\pi$
belonging to $U$ have no points in common, hence
\[
|P|\ge k(U):=\sum_{a\in U}k(a).
\]
This observation already proved the next lemma.

\begin{lemma}\label{independent}For each graph $G=(V,E)$ we have
\[
v(G)\ge\max\{k(U);\text{$U\subseteq V$ is independent}\},
\]
in particular $2\cdot\alpha(G)\leq v(G).$ 
\end{lemma}\qed

The numbers $k(a)$ can be used to show another lower
bound for the linear intersection number. Given any graph
$G=(V,E)$, define
\[
f(G):=\sum_{a\in V}k(a).
\]
Continuing our interpretation of $k(a)$, the number $f(G)$ is
something like a minimal number of flags required to realize
$G$ as an intersection graph.

\begin{lemma}If $G$ is a graph with $m$ edges, then
\[
v(G)\ge\frac{f(G)}{\omega(G)}\ge\frac{2m}{\omega(G)^2}.
\]
\end{lemma}

\textit{Proof.} Write $G=(V,E)$ and let $\pi=(P,V)$ be a linear
hypergraph with $G_{\pi}=G$ and $|P|=v(G)$.
Given any point $x\in P$, the line pencil $V_x$ is a clique
in $G$, hence $|V_x|\le\omega(G)$. We compute
\[
f(G)=\sum_{a\in V}k(a)\le\sum_{a\in V}|a|
=\sum_{x\in P}|V_x|\le|P|\omega(G)=v(G)\omega(G).
\]
For any graph $H$,
the vertices of $H$ may be partitioned in $\theta(H)$ cliques,
each of size at most $\omega(H)$, hence $\theta(H)\omega(H)$
is at least the number of vertices of $H$. Consequently,
\[
k(a)\ge\theta(G_a)\ge\frac{|G_a|}{\omega(G_a)}\ge
\frac{|G_a|}{\omega(G)}
\]
for any vertex $a\in V$. This implies
\[
f(G)=\sum_{a\in V}k(a)\ge\frac{1}{\omega(G)}\sum_{a\in V}|G_a|
=\frac{2m}{\omega(G)}.
\]
\qed

Consequently, the Erd\"os, Faber, Lov\'asz conjecture will be
true for graphs with $f(G)/\omega(G)\ge\chi(G)$. Since $\omega(G)
\le\chi(G)$, we obtain the next corollary.

\begin{corollary}Let $G$ be a graph with $\chi(G)\le\sqrt{f(G)}$.
Then $v(G)\ge\chi(G)$.\qed\end{corollary}

These little observations already indicate that it will be pretty
difficult to actually compute the number $v(G)$ for an arbitrary
graph $G$. For example, consider the complete $s$-partite graph
$K_{s,n}$ on $s\ge2$ sets each of size $n\ge2$. If $a$ is any vertex of
$K_{s,n}$, then $G_a=K_{s-1,n}$ and $k(a)=\theta(K_{s-1,n})=n$,
and Corollary \ref{independent} implies $v(K_{s,n})\ge n^2$.
Now, assume $v(K_{s,n})=n^2$, then it is possible to realize $K_{s,n}$
as the intersection graph of a linear hypergraph on $n^2$ points.
Each of the $s$ independent sets in $K_{s,n}$ will be a full parallel
pencil in this linear hypergraph, hence the linear hypergraph looks like
$s$ parallel pencils of an affine plane. Such linear hypergraphs are
sometimes called nets, and it is known that the existence of such
a net implies $s\le n+1$ and $s=n+1$ is possible if and only if there
exists a projective plane of order $n$.

In particular, an algorithm which computes $v(G)$ for each graph
$G$ is able to decide the existence of a projective plane of order $n$.
This fact pretty much implies that there will not exist an efficient
algorithm to compute $v(G)$ for an arbitrary graph $G$.

Of course, each graph is a subgraph of a certain complete graph. Therefore, it
is interesting to ask what happens if we remove an edge from a given graph.

\begin{lemma} Let $G$ be a graph and $G^-$ be some graph
obtained by deleting one
edge of $G$. Then $$v(G^-)\geq v(G)-1.$$
\end{lemma}

{\it Proof.} Let $x$ and $y$ be vertices of $G$ such that $e=\{x,y\}$
is an edge
of $G$, and let $G^-$ be the graph obtained from $G$ by deleting $e$.
Moreover, let
$\pi=(P,\L)$ be a linear hypergraph with $|P|=v(G^-)$
and $G_\pi=G^-.$ Then $x$ and $y$ are disjoint lines of $\pi.$ We introduce one
new point, denoted $\infty$, and define 
\[ \pi'=(P\cup\{\infty\}, ({\cal L}\backslash \{x,y\})\cup
\{x\cup\{\infty\},y\cup\{\infty\}\}).\]
Since the lines $  x\cup\{\infty\}$ and $y\cup\{\infty\}$ intersect, we obtain
$G=G_{\pi'}$, in particular $v(G)\leq |P\cup\{\infty\}|=v(G^-)+1.$\qed

\end{document}